\documentclass[11pt]{article}
\usepackage{graphicx}
\usepackage{amsmath}
\usepackage{amssymb}
\usepackage{amsfonts}
\usepackage{amsthm}
\usepackage{xcolor}
\usepackage{enumitem}
\usepackage{placeins}

\numberwithin{equation}{section}
\def\R{\mathbb{R}}

\def\neweq#1{\begin{equation}\label{#1}}
\def\endeq{\end{equation}}
\def\eq#1{(\ref{#1})}

\def\ee{{\bf e}}

\usepackage{hyperref}
\hypersetup{
    colorlinks=true,
    linkcolor=blue,
    filecolor=magenta,
    urlcolor=cyan,
}

\newtheorem{theorem}{Theorem}[section]

\newtheorem{proposition}[theorem]{Proposition}

\title{Some bounds for the eigenfunctions of the Stokes problem\\
under Navier boundary conditions in a cube}
\author{Gianni ARIOLI - Alessio FALOCCHI - Filippo GAZZOLA\\
{\small Dipartimento di Matematica\footnote{Dipartimento di Eccellenza MUR 2023-2027 (Italy)} \ -- Politecnico di Milano, Italy}}
\date{}
\begin{document}
\maketitle
\section{Introduction}
In the cube $\Omega=(0,\pi)^3$ we consider the following Stokes eigenvalue problem, studied in \cite{falgaz3}:
\begin{equation}\label{stokesphi}
	\left\{
	\begin{array}{ll}
		-\Delta \varphi=\lambda\varphi\quad &\mbox{in }\Omega\,,\\
		\nabla\cdot \varphi=0\quad &\mbox{in }\Omega\,,\\
		\varphi_{1}=\partial_x\varphi_{2}=\partial_x\varphi_{3}=0\quad &\mbox{on }\{0,\pi\}\times(0,\pi)\times(0,\pi)\,,\\
		\varphi_{2}=\partial_y\varphi_{1}=\partial_y\varphi_{3}=0\quad &\mbox{on }(0,\pi)\times\{0,\pi\}\times(0,\pi)\,,\\
		\varphi_{3}=\partial_z\varphi_{1}=\partial_z\varphi_{2}=0\quad &\mbox{on }(0,\pi)\times(0,\pi)\times\{0,\pi\}\,.
	\end{array}\right.
\end{equation}
We are interested in some properties of its eigenfunctions, that are computed explicitly in \cite{falgaz3}.

\begin{proposition}\label{autofunzioni}\cite{falgaz3}
	All the eigenvalues of \eqref{stokesphi} have finite multiplicity and can be ordered in an increasing divergent sequence
	$\{\lambda_k\}_{k\in \mathbb{N}_+}$, in which the eigenvalues are repeated according to their multiplicity.\par
	For $m,n,p\in \mathbb{N}_+$, the eigenfunctions
	\begin{equation}\label{X0}
		X_{0,n,p}(\xi):=\frac{2}{\sqrt{\pi^3(n^2+p^2)}}
		\left(\begin{array}{c}
			0\\
			p\sin(ny)\cos(pz)\\
			-n\cos(ny)\sin(pz)
		\end{array}\right)\ ,
	\end{equation}
	\begin{equation}\label{Y0}
		Y_{m,0,p}(\xi):=\frac{2}{\sqrt{\pi^3(m^2+p^2)}}
		\left(\begin{array}{c}
			-p\sin(mx)\cos(pz)\\
			0\\
			m\cos(mx)\sin(pz)
		\end{array}\right)\ ,
	\end{equation}
	\begin{equation}\label{Z0}
		Z_{m,n,0}(\xi):=\frac{2}{\sqrt{\pi^3(m^2+n^2)}}
		\left(\begin{array}{c}
			n\sin(mx)\cos(ny)\\
			-m\cos(mx)\sin(ny)\\
			0
		\end{array}\right)\ ,
	\end{equation}
	\begin{equation}\label{V}
		V_{m,n,p}(\xi):=\frac{2\sqrt2 \, \cos(mx)}{\sqrt{\pi^3(n^2+p^2)}}
		\left(\begin{array}{c}
			0\\
			p\sin(ny)\cos(pz)\\
			-n\cos(ny)\sin(pz)
		\end{array}\right)\ ,
	\end{equation}
	\begin{equation}\label{W}
		W_{m,n,p}(\xi):=\frac{2\sqrt2 }{\sqrt{\pi^3(m^2+n^2+p^2)(n^2+p^2)}}
		\left(\begin{array}{c}
			(n^2+p^2)\sin(mx)\cos(ny)\cos(pz)\\
			-mn\cos(mx)\sin(ny)\cos(pz)\\
			-mp\cos(mx)\cos(ny)\sin(pz)
		\end{array}\right)\ ,
	\end{equation}
	form an orthogonal basis in $U$ which is orthonormal in $H$. 
\end{proposition}
The statements proved in this paper will be used in \cite{arfalgaz}. We use the notation $\xi:=(x,y,z)\in\overline{\Omega}$.

\section{Bounds for the eigenfunctions}

We first compute the squared $L^\infty(\Omega)$-norm of the eigenfunctions in Proposition \ref{autofunzioni}.
	\begin{proposition}\label{proposition0}
	Let $m,n,p\in\mathbb{N}_+$, then
\begin{equation*}
	\begin{split}
		&\|X_{0,n,p}\|_{L^\infty(\Omega)}^2
		=\dfrac{4 \max\{n^2,p^2\}}{\pi^3(n^2+p^2)},\quad
		\|Y_{m,0,p}\|_{L^\infty(\Omega)}^2
		=\dfrac{4 \max\{m^2,p^2\}}{\pi^3(m^2+p^2)},\quad\|Z_{m,n,0}\|_{L^\infty(\Omega)}^2
		=\dfrac{4 \max\{m^2,n^2\}}{\pi^3(m^2+n^2)},\\
		&\|V_{m,n,p}\|_{L^\infty(\Omega)}^2=\dfrac{8 \max\{n^2,p^2\}}{\pi^3(n^2+p^2)},\quad
		\|W_{m,n,p}\|_{L^\infty(\Omega)}^2=\dfrac{8 \max\{(n^2+p^2)^2, m^2n^2,m^2p^2\}}{\pi^3(n^2+p^2)(m^2+n^2+p^2)}.
	\end{split}
\end{equation*}
\end{proposition}

Let $\ee=(a,b,c)\in\R^3$ be a unit vector, we now determine the maximum of the projection of the eigenfunctions in Proposition \ref{autofunzioni} onto $\ee$.

	\begin{proposition}\label{proposition1}
	Let $m,n,p\in\mathbb{N}_+$, $\ee=(a,b,c)\in\R^3$ be such that $a^2+b^2+c^2=1$. Then
	\begin{equation}\label{ts1}
		\begin{split}
			&\|\ee\cdot X_{0,n,p}\|_{L^\infty(\Omega)}^2
			=\dfrac{4 \max\{b^2p^2, c^2n^2\}}{\pi^3(n^2+p^2)},\\
			&\|\ee\cdot Y_{m,0,p}\|_{L^\infty(\Omega)}^2
			=\dfrac{4 \max\{a^2p^2,c^2m^2\}}{\pi^3(m^2+p^2)},\\
			&\|\ee\cdot Z_{m,n,0}\|_{L^\infty(\Omega)}^2
			=\dfrac{4 \max\{a^2n^2,b^2m^2\}}{\pi^3(m^2+n^2)},\\
			&\|\ee\cdot V_{m,n,p}\|_{L^\infty(\Omega)}^2=\dfrac{8 \max\{b^2p^2,c^2n^2\}}{\pi^3(n^2+p^2)},\\
			&\|\ee\cdot W_{m,n,p}\|_{L^\infty(\Omega)}^2=\dfrac{8 \max\{a^2(n^2+p^2)^2, b^2m^2n^2,c^2m^2p^2\}}{\pi^3(n^2+p^2)(m^2+n^2+p^2)}.
		\end{split}
	\end{equation}
\end{proposition}

In particular, \eqref{ts1}$_1$  is equivalent to
	\begin{equation*}
		\begin{split}
		&\max_{(y,z)\in[0,\pi]^2}\big[bp\sin (ny)\cos(pz)-cn\cos(ny)\sin(pz)\big]^2=\max\{b^2p^2,c^2n^2\},
		\end{split}
	\end{equation*}
whereas \eqref{ts1}$_5$ is equivalent to
	\begin{equation*}
	\begin{split}
		&\displaystyle\max_{(x,y,z)\in\overline\Omega}
		\big[a(n^2\!+\!p^2)\sin(mx)\cos(n y)\cos(pz)\!-\!bmn\cos(mx)\sin (ny)\cos(pz)\!-\!cmp\cos(mx)\cos(ny)\sin(pz)\big]^2\\
		&=\max\{a^2(n^2+p^2)^2,b^2m^2n^2,c^2m^2p^2\}\,.
	\end{split}
\end{equation*}
Hence, the proof of \eqref{ts1}$_{1,2,3,4}$ is almost immediate, while for \eqref{ts1}$_5$ it is lengthy and tedious.

Next, we bound sums of the squared projections  of the eigenfunctions along $\ee$. For instance, fix $\xi\in\overline{\Omega}$, $\mu>0$, and consider
\begin{eqnarray*}
\sum_{i=1}^h \dfrac{\big(\ee\cdot X_{0,n_i,p_i}(\xi)\big)^2}{(\lambda_i+\mu)^2} &=& \sum_{i=1}^h \dfrac{4}{\pi^3(n_i^2+p_i^2)}\dfrac{\big(bp_i\sin (n_iy)\cos(p_iz)-cn_i\cos(n_iy)\sin(p_iz)\big)^2}{(n_i^2+p_i^2+\mu)^2}\\
&\leq&\dfrac{4}{\pi^3}\sum_{i=1}^h \dfrac{\max\{b^2p_i^2,c^2n_i^2\}}{(n_i^2+p_i^2)(n_i^2+p_i^2+\mu)^2}
\le \dfrac{4}{\pi^3}\int_0^\infty\int_0^\infty \dfrac{\max\{b^2y^2,c^2x^2\}}{(x^2+y^2)(x^2+y^2+\mu)^2}dxdy,
\end{eqnarray*}
where we applied Proposition \ref{proposition1}. Similarly, for the sums related to the other four families of eigenfunctions in Proposition \ref{autofunzioni}.\par
Then we compute the integral bounding the sum. To do this, with the convention that $\arctan(+\infty)=\tfrac\pi2$, for $(b,c)\in\R^2$ such that $b^2+c^2\le1$ we introduce the function
\begin{eqnarray}\label{upsilon}
	\Upsilon(b,c) &:=&b^2\arctan\Big|\frac{c}{b}\Big| +c^2\arctan\Big|\frac{b}{c}\Big| +|bc|,
\end{eqnarray}
and we state the next result.
\begin{proposition}\label{proposition3}
	Let $\mu>0$, $m,n,p\in\mathbb{N}_+$ and $(a,b,c)\in\R^3$ be such that $a^2+b^2+c^2=1$. Then
	\begin{equation}\label{proposition4}
		\begin{split}
			&\int_0^\infty\int_0^\infty \dfrac{\max\{b^2y^2,c^2x^2\}}{(x^2+y^2)(x^2+y^2+\mu)^2}dxdy=\dfrac{1}{4\mu}\, \Upsilon(b,c)\\
			&\int_0^\infty\int_0^\infty\int_0^\infty \dfrac{\max\{b^2y^2,c^2x^2\}}{(x^2+y^2)(x^2+y^2+z^2+\mu)^2}dxdydz=
			\dfrac{\pi}{8\sqrt{\mu}}\, \Upsilon(b,c)\\
			&\int_0^\infty\int_0^\infty\int_0^\infty \dfrac{\max\{a^2(y^2+z^2)^2,b^2x^2y^2,c^2x^2z^2\}}{(y^2+z^2)(x^2+y^2+z^2)(x^2+y^2+z^2+\mu)^2}\, dxdydz\leq\frac{\pi}{12\sqrt{\mu}}\bigg[\pi\, a^2+ \dfrac{\Upsilon(b,c)}{2}\bigg].
		\end{split}
	\end{equation}
\end{proposition}

Our last computation deals with a suitable combination of the integrals in Proposition \ref{proposition3} in order to maximize it with respect to the components of the unit vector $\ee$.

\begin{proposition}\label{proposition44}
Let $k_0:=1+4\sqrt{2}\frac{\pi}{3}$, $k_1:=\frac{2}{3}\sqrt{2}\pi^2$ and
\begin{equation*}
\Gamma(a,b,c):=k_0\Upsilon(b,c)+\Upsilon(a,c)+\Upsilon(a,b)+k_1a^2,
\end{equation*}
where the function $\Upsilon(\cdot,\cdot)$ is defined in \eqref{upsilon}. Then it holds
\neweq{claimGamma}
\max_{a^2+b^2+c^2=1}\ \Gamma(a,b,c)\approx10.91\, .
\endeq
\end{proposition}

\section{Proofs of the bounds}\label{proofproposition}
\begin{proof}[\textbf{Proof of Proposition \ref{proposition0}.}] We maximize
	$$
	\xi\mapsto|X_{0,n,p}(\xi)|^2=\dfrac{4}{\pi^3(n^2+p^2)}[p^2\sin^2(ny)\cos^2(pz)+n^2\cos^2(ny)\sin^2(pz)]\qquad\xi\in\overline{\Omega}\, .
	$$
	By setting $Y=\sin^2(ny)$ and $Z=\sin^2(pz)$, this is equivalent to maximize
	$$\Gamma(Y,Z):=p^2Y(1-Z)+n^2(1-Y)Z$$
	over the square $(Y,Z)\in[0,1]^2$. The function $\Gamma$ is harmonic and strictly positive in $(0,1)^2$: therefore, it achieves the
	(positive) maximum on the boundary of the square $[0,1]^2$. Hence,
	\neweq{supX}
	\max_{(Y,Z)\in[0,1]^2}\ \Gamma(Y,Z)=\max\{n^2,p^2\},
	\endeq
	giving the first norm in Proposition \ref{proposition0}.
	We proceed similarly for the eigenfunctions $Y_{m,0,p}(\xi)$ and $Z_{m,n,0}(\xi)$.
	
	On the other hand, if we allow nonconstant coefficients in the linear combinations, we find the following relations between eigenfunctions:
	\neweq{relazioneV}
	V_{m,n,p}(\xi)=\sqrt2 \, \cos(mx)\, X_{0,n,p}(\xi)\qquad\forall\xi\in\Omega\, ,
	\endeq
	$$
	W_{m,n,p}(\xi)=\tfrac{n\sqrt2 \sqrt{m^2+n^2}}{\sqrt{(m^2+n^2+p^2)(n^2+p^2)}} \, \cos(pz)\, Z_{m,n,0}(\xi)
	-\tfrac{p\sqrt2 \sqrt{m^2+p^2}}{\sqrt{(m^2+n^2+p^2)(n^2+p^2)}} \, \cos(ny)\, Y_{m,0,p}(\xi)\quad\forall\xi\in\Omega\, .
	$$
	Then, for the eigenfunctions \eq{V}, from \eq{relazioneV} we infer that
	$$
	\|V_{m,n,p}\|_{L^\infty(\Omega)}^2=2\|X_{0,n,p}\|_{L^\infty(\Omega)}^2=\dfrac{8 \max\{n^2,p^2\}}{\pi^3(n^2+p^2)}.
	$$
	
	Finally, for the eigenfunctions \eq{W} we notice that
$$|W_{m,n,p}(\xi)|^2=\tfrac{8[(n^2+p^2)^2\sin^2(mx)\cos^2(ny)\cos^2(pz)+m^2n^2\cos^2(mx)\sin^2(ny)\cos^2(pz)+m^2p^2\cos^2(mx)\cos^2(ny)\sin^2(pz)]}
{\pi^3(n^2+p^2)(m^2+n^2+p^2)}
$$
and we have to maximize this function over $\overline{\Omega}$.
	By setting $X=\sin^2(mx)$, $Y=\sin^2(ny)$, $Z=\sin^2(pz)$, and
	$$
	\Gamma(\xi):=(n^2+p^2)^2X(1-Y)(1-Z)+m^2n^2(1-X)Y(1-Z)+m^2p^2(1-X)(1-Y)Z\, ,$$
	we are led to maximize $\Gamma$ over the cube $[0,1]^3$. Since $\Gamma$ is positive and harmonic, it achieves its maximum on the boundary of
	this cube and
	$$
	\max_{\xi\in[0,1]^3}\Gamma(\xi)=\max\{(n^2+p^2)^2,m^2n^2,m^2p^2\}\,
	$$
	completing the proof of Proposition \ref{proposition0}.
\end{proof}

\begin{proof}[\textbf{Proof of Proposition \ref{proposition1}.}] We begin with the proof of \eqref{ts1}$_1$. Let us introduce the function $M(y,z)=\big[bp\sin (ny)\cos(pz)-cn\cos(ny)\sin(pz)\big]^2$ with $b^2+c^2\leq 1$. The gradient of $M$ vanishes either
	when $M=0$ (and $M$ attains its minimum 0) or when $(y,z)\in[0,\pi]^2$ solve the homogeneous trigonometric system
	\begin{equation*}
		\begin{pmatrix}
			bp      & cn \\
			cn       & bp  \\
		\end{pmatrix}\begin{pmatrix}
			\cos(ny)\cos(pz)     \\
			\sin(ny)\sin(pz)      \\
		\end{pmatrix}=\begin{pmatrix}
			0    \\
			0     \\
		\end{pmatrix}\, .
	\end{equation*}
	If $(bp)^2\neq(cn)^2$ we find as possible extremal points
	$$
	\begin{cases}
		\cos(ny)\cos(pz) =0\\
		\sin(ny)\sin(pz) =0
	\end{cases}\qquad \Longrightarrow\qquad\cos(ny)=\sin(pz)=0\mbox{ or }\sin(ny)=\cos(pz)=0\, ,
	$$
	which proves the proposition. If $(bp)^2=(cn)^2$, then $M(y,z)=c^2n^2\sin^2(ny\mp pz)\leq c^2n^2$, proving again the proposition.
		The proof of \eqref{ts1}$_{2,3}$ follows similarly.
		
		To prove \eqref{ts1}$_4$ we recall \eqref{relazioneV}, so that
		$$
		(\ee\cdot V_{m,n,p}(\xi))^2=\frac{8 \cos^2(mx)}{\pi^3(n^2+p^2)}M(y,z);
		$$
		arguing on the function $M(y,z)$ as before we obtain  \eqref{ts1}$_4$.

The remaining part of the proof is devoted to \eqref{ts1}$_5$. We put $A:=a(n^2+p^2)$, $B:=-bmn$, $C:=-cmp$ and, over $\overline\Omega$, define the function
	\begin{equation*}
		K(\xi)=A\sin(mx)\cos(ny)\cos(pz)+B\cos(mx)\sin(ny)\cos(pz)+C\cos(mx)\cos(ny)\sin(pz)\, .
	\end{equation*}
	The proof follows if we show that
	\neweq{claimK}
	\|K\|_{L^\infty(\Omega)}=\max\{|A|,|B|,|C|\}\, .
	\endeq
	By the Weierstrass Theorem, $K$ achieves a maximum over $\overline\Omega$ and two different situations may occur.\par\noindent
	$\bullet$ {\bf The maximum is achieved on $\partial\Omega$.} Here, we notice that
	$$
	|K(0,y,z)|=|K(\pi,y,z)|=|B\sin(ny)\cos(pz)+C\cos(ny)\sin(pz)|\quad\forall(y,z)\in[0,\pi]^2\, ,
	$$
	$$
	|K(x,0,z)|=|K(x,\pi,z)|=|A\sin(mx)\cos(pz)+C\cos(mx)\sin(pz)|\quad\forall(x,z)\in[0,\pi]^2\, ,
	$$
	$$
	|K(x,y,0)|=|K(x,y,\pi)|=|A\sin(mx)\cos(ny)+B\cos(mx)\sin(ny)|\quad\forall(x,y)\in[0,\pi]^2\, .
	$$
	Then, by Proposition \ref{proposition1}, we know that
	$$
	\max_{(y,z)\in[0,\pi]^2}|K(0,y,z)|=\max_{(y,z)\in[0,\pi]^2}|K(\pi,y,z)|=\max\{|B|,|C|\}\, ,
	$$
	$$
	\max_{(x,z)\in[0,\pi]^2}|K(x,0,z)|=\max_{(x,z)\in[0,\pi]^2}|K(x,\pi,z)|=\max\{|A|,|C|\}\, ,
	$$
	$$
	\max_{(x,y)\in[0,\pi]^2}|K(x,y,0)|=\max_{(x,y)\in[0,\pi]^2}|K(x,y,\pi)|=\max\{|A|,|B|\}\, .
	$$
	Therefore, if the maximum of $|K|$ is achieved on $\partial\Omega$, then \eq{claimK} is proved.\par\noindent
	$\bullet$ {\bf The maximum is achieved in $\Omega$.} Let $\overline{\xi}\in\Omega$ be a maximum point for $|K(\xi)|$, then $\nabla K(\overline{\xi})=0$, i.e.,
	$$
	\left\{\begin{array}{lll}
		A\cos(mx)\cos(ny)\cos(pz)-B\sin(mx)\sin(ny)\cos(pz)-C\sin(mx)\cos(ny)\sin(pz)=0\, ,\\
		-A\sin(mx)\sin(ny)\cos(pz)+B\cos(mx)\cos(ny)\cos(pz)-C\cos(mx)\sin(ny)\sin(pz)=0\, ,\\
		-A\sin(mx)\cos(ny)\sin(pz)-B\cos(mx)\sin(ny)\sin(pz)+C\cos(mx)\cos(ny)\cos(pz)=0\, .
	\end{array}\right.
	$$
	Putting $X=mx$, $Y=ny$, $Z=pz$, the above system becomes
	\neweq{XYZ}
	\left\{\begin{array}{lll}
		A\cos(X)\cos(Y)\cos(Z)-B\sin(X)\sin(Y)\cos(Z)-C\sin(X)\cos(Y)\sin(Z)=0\, ,\\
		A\sin(X)\sin(Y)\cos(Z)-B\cos(X)\cos(Y)\cos(Z)+C\cos(X)\sin(Y)\sin(Z)=0\, ,\\
		A\sin(X)\cos(Y)\sin(Z)+B\cos(X)\sin(Y)\sin(Z)-C\cos(X)\cos(Y)\cos(Z)=0\
	\end{array}\right.
	\endeq
	and, from now on, the original parameters ($a,b,c,x,y,z,m,n,p$) will not appear any more.
	We compute the values of $K=K(X,Y,Z)$ for any triplet $(X,Y,Z)\in\R^3$ satisfying \eq{XYZ}.
	To this end, we distinguish several cases, with some subcases.\par
	\underline{\em Case 1:} $\cos(X)\cos(Y)\cos(Z)=0$.\par
	If $\cos(X)=0$ and if the critical point $\overline{\xi}\in\Omega$ is a maximum for $|K|$, then
	$|K(\overline{\xi})|=|A\cos(Y)\cos(Z)|\le|A|$, with equality for suitable $(Y,Z)$, so that \eq{claimK} is satisfied.
	Similarly, if $\cos(Y)=0$ we find that $|K(\overline{\xi})|\le|B|$, while if $\cos(Z)=0$ we find that $|K(\overline{\xi})|\le|C|$.
	Therefore, in Case 1, \eq{claimK} is satisfied.\par
	\underline{\em Case 2:} $\cos(X)\cos(Y)\cos(Z)\neq0$.\par
	After dividing by $\cos(X)\cos(Y)\cos(Z)$ each of the equations in \eq{XYZ}, we obtain the system
	\neweq{XYZ1}
	\left\{\begin{array}{lll}
		B\tan(X)\tan(Y)+C\tan(X)\tan(Z)=A\, ,\\
		A\tan(X)\tan(Y)+C\tan(Y)\tan(Z)=B\, ,\\
		A\tan(X)\tan(Z)+B\tan(Y)\tan(Z)=C\, .
	\end{array}\right.
	\endeq
	In order to solve \eq{XYZ1} we distinguish two subcases.\par
	\underline{\em Case 2.1:} $ABC=0$.\par
	If $A=0$, then \eq{XYZ1}$_2$-\eq{XYZ1}$_3$ yield $B^2=C^2>0$ so that $C=\pm B$ and $K$ may be written as
	$$K(\xi)=B\cos(X)[\sin(Y)\cos(Z)\pm\cos(Y)\sin(Z)]=B\cos(X)\sin(Y\pm Z)\, ,$$
	implying that $|K(\xi)|\le|B|$ for all $\xi\in\overline\Omega$ and showing that \eq{claimK} holds. Similarly, if $B=0$ then
	$\|K\|_{L^\infty(\Omega)}\le|C|$, while if $C=0$ then $\|K\|_{L^\infty(\Omega)}\le|A|$. Hence, if $ABC=0$, then \eq{claimK} holds.\par
	\underline{\em Case 2.2:} $ABC\neq0$.\par
	In this case, the coefficient matrix of the (linear) system \eq{XYZ1} is nonsingular and \eq{XYZ1} is uniquely solvable in terms of the
	three unknowns $\tan(X)\tan(Y)$, $\tan(X)\tan(Z)$, $\tan(Y)\tan(Z)$. We find
	\neweq{tre}
	\tan(X)\tan(Y)=\frac{A^2\!+\!B^2\!-\!C^2}{2AB},\ \tan(X)\tan(Z)=\frac{A^2\!+\!C^2\!-B^2}{2AC},\ \tan(Y)\tan(Z)=\frac{B^2\!+\!C^2\!-\!A^2}{2BC}.
	\endeq
	We now claim that
	$$
	\tan(X)\tan(Y)\tan(Z)=0\ \Longrightarrow\ \mbox{\eq{claimK} holds}.
	$$
	Indeed, if $\tan(X)=0$ from \eq{tre}$_1$-\eq{tre}$_2$ we infer that $A^2+B^2-C^2=A^2+C^2-B^2=0$, which yields $B^2=C^2$ and, at the corresponding
	critical point $\overline{\xi}$, we have that
	$$|K(\overline{\xi})|=|B|\, |\sin(Y)\cos(Z)\pm\cos(Y)\sin(Z)|=|B|\, |\sin(Y\pm Z)|\le|B|\, ,$$
	with equality for suitable $(Y,Z)$: hence, \eq{claimK} holds. We argue similarly when $\tan(Y)=0$ or $\tan(Z)=0$.\par
	Therefore, we may assume that $\tan(X)\tan(Y)\tan(Z)\neq0$ and, in such case, we find
	$$\tan^2(X)=\frac{(A^2+B^2-C^2)(A^2+C^2-B^2)}{2A^2(B^2+C^2-A^2)},\quad\tan^2(Y)=\frac{(A^2+B^2-C^2)(B^2+C^2-A^2)}{2B^2(A^2+C^2-B^2)},$$
	\neweq{tansq}
	\tan^2(Z)=\frac{(A^2+C^2-B^2)(B^2+C^2-A^2)}{2C^2(A^2+B^2-C^2)}.
	\endeq
	From \eq{tansq} we infer that, necessarily, $(A^2+C^2-B^2)(B^2+C^2-A^2)(A^2+B^2-C^2)>0$,
	namely
	\neweq{alternative}
	\mbox{either the three terms are all positive or two of them are negative and the third is positive.}
	\endeq
	From \eq{tansq} we also deduce that
	$$\cos^2(X)=\frac{2 A^2 \left(B^2\!+\!C^2\!-\!A^2\right)}{2(A^2B^2\!+\!A^2C^2\!+\!B^2C^2)\!-\!A^4\!-\!B^4\!-\!C^4},\quad
	\cos^2(Y)=\frac{2 B^2 (A^2\!+\!C^2\!-\!B^2)}{2(A^2B^2\!+\!A^2C^2\!+\!B^2C^2)\!-\!A^4\!-\!B^4\!-\!C^4},$$
	\neweq{cossq}
	\cos^2(Z)=\frac{2C^2(A^2\!+\!B^2\!-\!C^2)}{2(A^2B^2\!+\!A^2C^2\!+\!B^2C^2)\!-\!A^4\!-\!B^4\!-\!C^4}.
	\endeq
	Since the three denominators in \eq{cossq} coincide, together with \eq{alternative}, this enables us to conclude that
	\neweq{allpositive}
	A^2+C^2>B^2\, ,\quad B^2+C^2>A^2\, ,\quad A^2+B^2>C^2\, ,\quad 2(A^2B^2\!+\!A^2C^2\!+\!B^2C^2)>A^4+B^4+C^4\, .
	\endeq
	Assume that $A^2,B^2\le C^2$ (the other cases being similar), put $\alpha=A^2/C^2\le1$ and $\beta=B^2/C^2\le1$, so that $0<\alpha,\beta\le1$
	and the four conditions in \eq{allpositive} reduce to the third condition, namely we obtain the triangle $\Theta$ defined by
	\neweq{constraint}
	\Theta=\{(\alpha,\beta)\in(0,1]^2;\ \alpha+\beta>1\}\, .
	\endeq
	
	We then notice that $\|K\|_{L^\infty(\Omega)}$ is the squared root of the maximum of the function
	\begin{equation*}
		\begin{split}
			K^2(\xi)=&\cos^2(X)\cos^2(Y)\cos^2(Z)[A\tan(X)+B\tan(Y)+C\tan(Z)]^2\, \\
			=&\cos^2(X)\cos^2(Y)\cos^2(Z)\big[A^2\tan^2(X)+B^2\tan^2(Y)+C^2\tan^2(Z)+\\
			&+2AB\tan(X)\tan(Y)+2BC\tan(Y)\tan(Z)+2AC\tan(X)\tan(Z)\big].
		\end{split}
	\end{equation*}
	An interior critical point for $K$ is also an interior critical point for $K^2$. Therefore, if $ABC\neq0$ (the case $ABC=0$ has already been
	discussed) and if $K^2$ achieves its maximum at some interior point $(X,Y,Z)$, then \eq{tre}-\eq{tansq}-\eq{cossq} hold. Hence, in such point
	we have
	\begin{equation*}
		\begin{split}
			K^2(\xi)=&\tfrac{8A^2B^2C^2(B^2\!+\!C^2\!-\!A^2)(A^2\!+\!C^2\!-\!B^2)(A^2\!+\!B^2\!-\!C^2)}{[2(A^2B^2\!+\!A^2C^2\!+\!B^2C^2)\!-\!A^4\!-\!B^4\!-\!C^4]^3}\\
			&\times\Big[\tfrac{(A^2+B^2-C^2)^2(A^2+C^2-B^2)^2+(A^2+B^2-C^2)^2(B^2+C^2-A^2)^2+
				(A^2+C^2-B^2)^2(B^2+C^2-A^2)^2}{2(B^2+C^2-A^2)(A^2+C^2-B^2)(A^2+B^2-C^2)}+A^2\!+\!B^2\!+\!C^2\Big]\\
		\end{split}
	\end{equation*}
	
	With the previous change of variables and some lengthy computations, this expression simplifies and we may rewrite it as
	\begin{equation*}
		\frac{K^2(\xi)}{C^2}=D(\alpha,\beta):=\frac{4\alpha\beta}{2(\alpha\beta+\alpha+\beta)-\alpha^2-\beta^2-1}.
	\end{equation*}
	In the interior of the triangle $\Theta$ defined by \eq{constraint}, $\nabla D\neq0$ so that $D$ achieves its maximum over $\overline\Theta$
	on $\partial\Theta$ (in fact, for $\alpha+\beta=1$) and $\max_{\partial\Theta}D=1$. Therefore, \eq{claimK} holds also in this case.
\end{proof}
\noindent
\begin{proof}[\textbf{Proof of Proposition \ref{proposition3}.}]  Let us first compute the integral
	\begin{equation}
		\begin{split}
			&2\int_0^{\pi/2}\max\{b^2\cos^2\theta,c^2\sin^2\theta\}d\theta=
			2\bigg[b^2\int_0^{\arctan|\frac{c}{b}|}\cos^2\theta d\theta+c^2\int_{\arctan|\frac{c}{b}|}^{\frac{\pi}{2}}\sin^2\theta d\theta\bigg]=\\
			& b^2\int_0^{\arctan|\frac{c}{b}|}(1+\cos2\theta)d\theta
			+c^2\int_{\arctan|\frac{c}{b}|}^{\frac{\pi}{2}}(1-\cos2\theta) d\theta=\\
			&b^2\arctan\left|\frac{c}{b}\right| +\frac{b^2}{2}\sin\bigg(2\arctan\bigg|\frac{c}{b}\bigg|\bigg)
			+c^2\left(\frac\pi2 -\arctan\left|\frac{c}{b}\right|\right) -\frac{c^2}{2}\sin\bigg(2\arctan\bigg|\frac{c}{b}\bigg|\bigg)=\\
			& b^2\arctan\left|\frac{c}{b}\right| +c^2\left(\frac\pi2 -\arctan\left|\frac{c}{b}\right|\right) +|bc|=
			b^2\arctan\Big|\frac{c}{b}\Big| +c^2\arctan\Big|\frac{b}{c}\Big| +|bc|=\Upsilon(b,c).
		\end{split}
	\end{equation}

	To prove \eq{proposition4}$_1$ we switch to polar coordinates and we find
	\begin{equation*}
		\int_0^\infty\int_0^\infty \dfrac{\max\{b^2y^2,c^2x^2\}}{(x^2+y^2)(x^2+y^2+\mu)^2}dxdy=
		\frac{\Upsilon(b,c)}{2}\, \int_0^\infty\dfrac{\rho}{(\rho^2+\mu)^2}d\rho=\dfrac{1}{4\mu}\Upsilon(b,c).
	\end{equation*}
	
	To get \eq{proposition4}$_2$ we use \eq{proposition4}$_1$ as follows:
	$$
	\int_0^\infty\int_0^\infty\int_0^\infty \dfrac{\max\{b^2y^2,c^2x^2\}}{(x^2+y^2)(x^2+y^2+z^2+\mu)^2}dxdydz =\frac{\Upsilon(b,c)}{4}
	\int_0^\infty\dfrac{dz}{\mu+z^2}\, = \dfrac{\pi}{8\sqrt{\mu}}\Upsilon(b,c)
	$$
	because $\mu+z^2$ replaces $\mu$ in the inner $(x,y)$-integrals.
		
	To prove \eq{proposition4}$_3$ we need to bound
	$$
	I:=\int_0^\infty\int_0^\infty\int_0^\infty \dfrac{\max\{a^2(y^2+z^2)^2,b^2x^2y^2,c^2x^2z^2\}}{(y^2+z^2)(x^2+y^2+z^2)(x^2+y^2+z^2+\mu)^2}\, dxdydz.
	$$
	By putting $x=\rho t$, $y=\rho\cos\theta$, $z=\rho\sin\theta$, the Jacobian is $\rho^2$ and we obtain
	\begin{eqnarray*}
		I &=& \int_0^\infty\frac{dt}{t^2+1}\int_0^\infty\frac{\rho^2\, d\rho}{(\rho^2t^2+\rho^2+\mu)^2}
		\int_0^{\pi/2}\max\{a^2,b^2t^2\cos^2\theta,c^2t^2\sin^2\theta\}d\theta\\
		(s=\rho\sqrt{t^2+1}) &=& \int_0^\infty\frac{dt}{(t^2+1)^{5/2}}\int_0^\infty\frac{s^2\, ds}{(s^2+\mu)^2}
		\int_0^{\pi/2}\max\{a^2,b^2t^2\cos^2\theta,c^2t^2\sin^2\theta\}d\theta\\
		&\le& \left[\frac{\arctan\tfrac{s}{\sqrt{\mu}}}{2\sqrt{\mu}}-\frac{s}{2(s^2+\mu)}\right]_0^\infty
		\int_0^\infty\!\!\frac{dt}{(t^2+1)^{5/2}}\int_0^{\pi/2}\!\!\big[a^2+t^2\max\{b^2\cos^2\theta,c^2\sin^2\theta\}\big]d\theta\\
		&=& \frac{\pi}{4\sqrt{\mu}}\!\!\left[\frac{\pi\, a^2}{2}\int_0^\infty\!\!\frac{dt}{(t^2+1)^{5/2}}+
		\dfrac{\Upsilon(b,c)}{2}\int_0^\infty\frac{t^2\, dt}{(t^2+1)^{5/2}}\right]\\
		&=& \frac{\pi}{4\sqrt{\mu}}\left\{\frac{\pi\, a^2}{2}\left[\frac{3t+2t^3}{3(1+t^2)^{3/2}}\right]_0^\infty+
		\dfrac{\Upsilon(b,c)}{2}\left[\frac{t^3}{3(1+t^2)^{3/2}}\right]_0^\infty\right\}
		= \frac{\pi}{12\sqrt{\mu}}\bigg[\pi\, a^2+ \dfrac{\Upsilon(b,c)}{2}\bigg]\, ,
	\end{eqnarray*}
	and \eq{proposition4}$_3$ is proved.
\end{proof}
\noindent
\begin{proof}[\textbf{Proof of Proposition \ref{proposition44}.}]
Let $k_2:=\frac{\pi+2}{4}k_0$; we introduce the function
	$$
	G(s):=\frac{2(1-s^2)\arctan s+2s+k_1-k_2-\pi/2}{1+2s^2}\qquad(s\ge0)
	$$
	and we find its maximum. Hence, we compute the first derivative
	$$
	G'(s)=\frac{1-(s^3+s)(3\arctan s+k_1-k_2-\pi/2)-2s^4}{(1+2s^2)^2(1+s^2)}\qquad(s\ge0).
	$$
	It is possible to prove that $G'(s)$ admits a unique zero at $\sigma \approx0.672$, where $G(s)$ has a unique maximum. Then we put
	\neweq{G0}
	G_0:=\max_{s\ge0}G(s)=G(\sigma)\approx0.435\, .
	\endeq

	By possibly replacing $a,b,c$ with their absolute values, we may restrict our attention to the first octant where $a,b,c\ge0$.
	We first notice that
	\neweq{Upsilon1}
	\Upsilon(b,c)=\Upsilon(c,b)\, ,\qquad\Upsilon(b,b)=\left(1+\frac{\pi}{2}\right)b^2\, ,\qquad\Upsilon_b(b,c)=\frac{2c^3}{b^2+c^2}+2b
	\, \arctan\frac{c}{b}\ge0
	\endeq
	and, similarly, for $\Upsilon(a,b)$ and $\Upsilon(a,c)$. By the Lagrange multiplier method, in the maximum point there exists $\lambda\neq0$
	such that $\nabla\Gamma=\lambda\nabla g$, in particular that
	\begin{eqnarray*}
		k_0\left(\frac{c^3}{b^2+c^2}+b\, \arctan\frac{c}{b}\right)+\frac{a^3}{b^2+a^2}+b\, \arctan\frac{a}{b} &=& \lambda b\\
		k_0\left(\frac{b^3}{b^2+c^2}+c\, \arctan\frac{b}{c}\right)+\frac{a^3}{a^2+c^2}+c\, \arctan\frac{a}{c} &=& \lambda c\, .
	\end{eqnarray*}
	By \eq{Upsilon1} this implies that $c=b$ and the original problem \eq{claimGamma} becomes
	$$
	\mbox{maximize}\quad k_0\left(1+\frac{\pi}{2}\right)b^2+2\Upsilon(a,b)+k_1a^2\quad\mbox{under the constraint}\quad a^2+2b^2=1
	$$
	or, equivalently,
	$$
	\mbox{maximize}\qquad\Lambda(a):=k_0\left(1+\frac{\pi}{2}\right)\frac{1-a^2}{2}
	+2\Upsilon\left(a,\frac{\sqrt{1-a^2}}{\sqrt2 }\right)+k_1a^2\qquad\mbox{for }a\in[0,1]\, .
	$$
	After putting $\Lambda_0(a):=\Lambda(a)-k_2$, we are led to maximize
	$$
	a\mapsto2a^2\arctan\frac{\sqrt{1-a^2}}{a\sqrt2}+(1-a^2)\arctan\frac{a\sqrt2}{\sqrt{1-a^2}}+a\sqrt2\sqrt{1-a^2}+(k_1-k_2)a^2
	\qquad\mbox{for }a\in[0,1]\, .
	$$
	Next, we put
	$$
	s:=\frac{\sqrt{1-a^2}}{a\sqrt{2}}\ \Longrightarrow\ a^2=\frac{1}{1+2s^2}\, ,\quad1-a^2=\frac{2s^2}{1+2s^2}\, ,\quad
	a\sqrt{1-a^2}=\frac{s\sqrt2 }{1+2s^2}\, ,
	$$
	so that the target becomes to maximize, for $s\geq 0$,
	\begin{eqnarray*}
		F(s) &:=& \frac{2}{1+2s^2}\arctan s+\frac{2s^2}{1+2s^2}\arctan\frac{1}{s}+\frac{2s}{1+2s^2}+\frac{k_1-k_2}{1+2s^2}\\
		&=& \frac{2}{1+2s^2}\arctan s+\frac{2s^2}{1+2s^2}\left[\frac\pi2-\arctan s\right]+\frac{2s+k_1-k_2}{1+2s^2}\\
		&=& \frac{2(1-s^2)\arctan s+2s+k_1-k_2-\pi/2}{1+2s^2}+\frac\pi2=G(s)+\frac\pi2\, .
	\end{eqnarray*}
	
	By undoing the changes of variables, and recalling \eq{G0}, we find
	$$
	\max_{s\ge0}F(s)=G_0+\frac\pi2\ \Longrightarrow\ \max_{0\le a\le1}\Lambda_0(a)=G_0+\frac\pi2\ \Longrightarrow\
	\max_{0\le a\le1}\Lambda(a)=G_0+\frac\pi2 +\frac{2+\pi}4 \left(1+4\sqrt{2}\frac{\pi}{3}\right)\, ,
	$$
	which proves \eq{claimGamma}.
\end{proof}
\end{document}